\documentstyle{amsart}
\newtheorem{thm}{Theorem}[section]
\newtheorem{pro}[thm]{Proposition}
\newtheorem{lem}[thm]{Lemma}
\theoremstyle{definition}
\newtheorem{defn}{Definition}[section]
\numberwithin{equation}{section}
\def\dj{d\kern-.30em\raise1.25ex\vbox{\hrule width .3em height .03em}}
\def\Dj{D\kern-.75em\raise0.75ex\vbox{\hrule width .3em height .03em}
\kern.03em}
\def\sstar{{\raise0.2ex\hbox{$\scriptstyle\star$}}}
\newsymbol\restriction 1316
\newsymbol\between 1347
\def\restr{{\restriction}}
\newcommand{\hor}{\mbox{\family{euf}\shape{n}\selectfont hor}}
\newcommand{\inv}{i\!\hspace{0.8pt}n\!\hspace{0.6pt}v}
\newcommand{\pre}{\pi}
\newcommand{\nat}{\widehat{\cal{R}}}
\newcommand{\grten}{\mathbin{\widehat{\otimes}}}
\newcommand{\vh}{\mbox{\family{euf}\shape{n}\selectfont vh}}
\newcommand{\ad}{\mbox{\shape{n}\selectfont ad}}
\newcommand{\id}{\mbox{\shape{n}\selectfont id}}
\newcommand{\Ad}{\varpi}
\def\Uni{\between}
\def\bla#1{$(${\it #1\/{}}$)$}
\newcommand{\lie}{\mbox{\shape{n}\selectfont lie}}
\newcommand{\kr}{\ker}
\newcommand{\e}{\epsilon}
\newcommand{\I}{\Upsilon}
\newcommand{\Sum}{\displaystyle{\sum}}
\newcommand{\Bcap}{\textstyle\bigcap}
\begin{document}
\title[Differential Structures]{ON DIFFERENTIAL
STRUCTURES ON  QUANTUM \\ PRINCIPAL BUNDLES}
\author{Mi\'co {\Dj}ur{\Dj}evi\'c}
\address{Instituto   de   Matematicas,   UNAM,    Area    de    la
Investigacion Cientifica, Circuito Exterior, M\'exico DF, CP
04510, M\'EXICO\newline
\indent {\it Written In}\newline
\indent Centro de Investigaciones Teoricas, UNAM, Facultad de Estudios
Superiores Cuautitlan, Cuautitlan Izcalli, M\'EXICO}
\maketitle
\begin{abstract}
   A constructive approach to differential calculus on
quantum principal bundles is presented. The calculus on the bundle
is built in  an  intrinsic  manner,  starting  from  given  graded
(differential) *-algebras representing  horizontal  forms  on  the
bundle and differential forms on the base manifold, together with
a family of antiderivations acting on  horizontal  forms,  playing
the role of covariant derivatives of regular connections. In  this
conceptual framework,  a  natural  differential  calculus  on  the
structure quantum group is described.
\end{abstract}
\tableofcontents
\section{Introduction}
     The  aim  of  this  study  is  to  present   some   algebraic
constructions  related  to  differential   calculus   on   quantum
principal bundles.

     The paper is logically based on a general theory  of  quantum
principal bundles, developed in \cite{D1,D2}.

     As far as quantum principal  bundles  over  classical  smooth
manifolds are concerned, it is possible to construct, in a natural
manner, a differential calculus on them. The algebra
of  differential  forms  on  the  bundle  can  be  constructed  by
combining standard differential forms on the  base  manifold  with
appropriate differential calculus on the structure quantum group,
such  that  every  local  trivialization  of  the  bundle  can  be
``extended'' to a local trivialization  of  the corresponding
calculus.

{\renewcommand{\thepage}{}
However,
this  local  triviality   property   implies   relatively   strong
constraints for the algebra of differential
forms on the structure quantum group. At the first-order level, in
the class  of  left-covariant   differential  structures  (on  the
structure quantum group), there exists the unique minimal  element
satisfying mentioned constraints. This  (first-order) calculus is
also *-covariant and bicovariant. If the  higher-order
differential calculus on the structure group is described  by  the
corresponding   universal   envelope,   then   all   compatibility
conditions are resolved  already  at  the  first-order  level. The
same situation holds if the higher-order calculus is described  by
the corresponding exterior algebra \cite{W-diff}.

In summary,   quantum  principal  bundles   over   classical
smooth manifolds are sufficiently ``structuralized'' geometrical
objects.
This opens the possibility to construct
the whole differential calculus in an intrinsic manner,
starting from the idea of local triviality.

     On the other hand, in the theory of general quantum principal
bundles (over  quantum  spaces)  we  meet  just  the   oposite
situation.
{}From the conceptual point of view, the most  natural  approach  to
differential calculus is to start  from  appropriate  differential
*-algebras representing differential calculi on the bundle and
the structure quantum group.  Then  quantum  counterparts  of  all
basic entities appearing in the classical formalism can be derived
from these algebras. In particular, differential forms on the base
manifold can be  viewed  as  differential  forms  on  the  bundle,
invariant under the ``pull back'' induced by the right action map.

     In  this  paper  a  constructive  approach  to   differential
calculus on general quantum principal bundles will  be  presented.
This approach incorporates  some  ideas  of \cite{D1} into
the  general
quantum context. Technically everything can  be  considered  as  a
variation  of a theme presented in \cite{D2}--Subsection 6.6.
The
notation and terminology introduced in \cite{D2} will be  followed
here.

Let $G$ be a compact matrix quantum group \cite{W-cmpg}. Let $\cal
A$ be the Hopf *-algebra of polynomial functions on $G$. We  shall
denote by $\phi\colon\cal A\rightarrow\cal A\otimes\cal A$,
$\e\colon\cal A\rightarrow \Bbb{C}$ and
$\kappa\colon\cal A\rightarrow\cal A$
the coproduct, counit and the antipode respectively.
}

Let us  assume  that  a  complete  differential  calculus  on  the
structure quantum group $G$ is
specified    by    the    universal     differential      envelope
$\Gamma^{\wedge}$   of  a
first-order bicovariant *-calculus $\Gamma$.

   Let $P=(\cal{B},i,F)$ be a quantum principal $G$-bundle
over a quantum
space $M$, which is formally represented by a *-algebra $\cal{V}$. Here,
$\cal{B}$ formally represents $P$ while
$i\colon \cal{V}\rightarrow\cal{B}$ is the dualized projection of $P$ on
$M$   and   $F\colon \cal{B}\rightarrow\cal{B}\otimes\cal{A}$   is   the
dualized right action of $G$ on $P$.
Let $\Omega(P)$ be a graded-differential *-algebra
representing
the  complete  differential  calculus on  the  bundle.   Let
$\hor(P)\subseteq\Omega(P)$  be the corresponding
*-subalgebra of horizontal forms.

Every connection $\omega\colon \Gamma_{\inv}\rightarrow\Omega(P)$
(where
$\Gamma_{\inv}$ is the space of left-invariant elements of
$\Gamma$) on $P$  is  completely  determined  by  the
corresponding  covariant  derivative
$D_{\omega}\colon \hor(P)\rightarrow \hor(P)$.
This map is right-covariant and
$(D_{\omega}\restr\Omega(M))=d_M,$ where
$\Omega(M)\subseteq\Omega(P)$ is the  graded-differential *-algebra
representing  differential  forms  on  the  base  space $M$,   and
$d_M\colon \Omega(M)\rightarrow\Omega(M)$
is     the      corresponding
differential. If $\omega$ is regular then
$D_{\omega}$ is hermitian  and  satisfies  the
graded  Leibniz  rule.  The
structure of $\Omega(P)$ is completely encoded in $\bigl\{
\hor(P),
F^\wedge,D_\omega,\Gamma_{\inv}^\wedge\bigr\}$ (where
$\Gamma_{\inv}^{\wedge}\subseteq\Gamma^\wedge$ is
a differential *-subalgebra consisting of left-invariant
elements).

     The starting point for considerations of this paper  consists
of a graded *-algebra $\hor_P$ (the elements of which are
interpretable as horizontal  forms
on $P$), a *-homo\-mor\-phism
$F^{\sstar}\colon \hor_P\rightarrow \hor_P\otimes\cal{A}$ (playing the
role of the induced action of $G$ on horizontal forms)
and  the  graded  subalgebra $\Omega_M\subseteq \hor_P$ consisting
of $F^{\sstar}$-invariant elements (representing
differential forms on $M$)     endowed      with      a      differential
$d_M\colon \Omega_M\rightarrow\Omega_M$.
 Then  counterparts   of   covariant
derivatives (of regular connections) can be defined
as operators acting in $\hor_P$, and possessing  the  above  mentioned
characteristic  properties.  Such   operators   will   be   called
{\it preconnections}.

     Section 2 is devoted to the study of  elementary  properties
of preconnections.
     Starting from the space of preconnections, it is possible  to
construct, in a natural manner,  a  graded-differential *-algebra
$\Omega_P$ imaginable as  consisting  of
differential  forms  on $P$,
together with an appropriate  bicovariant  first-order
*-calculus
$\Psi$ on $G$. These constructions will be
presented in  Section 3. The constructed algebra $\Omega_P$
contains $\hor_P$ as its graded *-subalgebra. Actually $\hor_P$
coincides with the graded *-subalgebra
representing horizontal forms in the
general theory. The map
$F^{\sstar}\colon \hor_P\rightarrow \hor_P\otimes\cal{A}$
coincides with the associated action
of  $G$  on  horizontal  forms
(defined in the framework of the general theory). This
implies that $\Omega_M$ consists precisely of
$\widehat{F}$-invariant forms,
where $\widehat{F}\colon\Omega_P
\rightarrow\Omega_P\grten\Psi^\wedge$
is the natural
graded-differential extension of $F$ (the ``pull back'' map).
Further, regular connections on $P$
(relative  to $\bigl\{\Omega_P,\Psi^\wedge\bigr\}$)  are
multiplicative  and  there
exists a natural correspondence between  regular  connections  and
preconnections on $P$ (interpretable  as labeling of regular
connections  by   the   corresponding   operators   of   covariant
derivative). In such a way, a circle will be closed.

    Finally, in Section 4 some concluding remarks are made.

     The paper  ends  with  a  technical  appendix  in  which  the
construction of the calculus on the bundle  is  sketched in  the
case when the higher-order differential calculus on the  structure
group is  described  by  the  corresponding  bicovariant  exterior
algebra \cite{W-diff}. The  calculus  can  be
constructed essentially in  the  same  way  as  in  the
``universal envelope'' case.
\section{Preconnections and Their Elementary Properties}

Let $M$ be a quantum space, represented by a *-algebra $\cal{V}$.
Let $G$  be  a  compact  matrix  quantum   group \cite{W-cmpg} and
let $P=(\cal{B},i,F)$ be
a quantum principal $G$-bundle \cite{D2}
over $M$.  The action
$F\colon\cal{B}\rightarrow\cal{B}\otimes
\cal{A}$ is free, in the
sense  that  for
each $a\in\cal{A}$ there exist elements $b_i,q_i\in\cal{B}$
such that
\begin{equation}\label{free}
\sum_iq_iF(b_i)=1\otimes a.
\end{equation}
     Let
$$\hor_P=\sideset{}{^\oplus}\sum_{k\geq 0} \hor_P^k$$
be a graded *-algebra such that
$\cal{B}=\hor_P^0$ and let
$F^{\sstar}\colon \hor_P\rightarrow \hor_P\otimes\cal{A}$
be a grade-preserving *-homomorphism extending $F$ and satisfying
\begin{align}
(\id\otimes\phi)F^{\sstar}
&=(F^{\sstar}\otimes \id)F^{\sstar}\label{coas}\\
(\id\otimes \e)F^{\sstar}&=\id. \label{ide}
\end{align}
Geometrically, $F^{\sstar}$  determines a  (left)  action  of  $G$
by
``automorphisms'' of $\hor_P$.

The elements of $\hor_P$  will be  interpreted
as horizontal differential forms on $P$.
Let $\Omega_M\subseteq \hor_P$ be the graded *-subalgebra
consisting of $F^{\sstar}$-invariant  elements.  The  elements  of
$\Omega_M$ will be interpreted as  differential  forms  on  the  base
manifold $M$
(these interpretations will be completely justified
after constructing the calculus on the bundle $P$).

\begin{lem}\label{lem:21} Let us assume that a linear map
$\Delta\colon \hor_P\rightarrow \hor_P$ is given such that
\begin{gather}
F^{\sstar}\Delta=(\Delta\otimes \id)F^{\sstar}\label{FDel}\\
\Delta(\Omega_M)=\{0\}.\label{DelW}
\end{gather}

     \bla{i} If $\Delta$ is an (odd) antiderivation then there exists
the unique $\varrho\colon \cal{A}\rightarrow \hor_P$
such that
\begin{equation}\label{strD1}
\Delta(\varphi)=(-1)^{\partial\varphi}\sum_k\varphi_k\varrho(c_k)
\end{equation}
for each $\varphi\in \hor_P$,  where $\Sum_k\varphi_k\otimes c_k
=F^{\sstar}(\varphi)$.
The following identity holds
\begin{equation}\label{com1}
\varrho(a)\varphi=(-1)^{\partial\varphi}\sum_k\varphi_k\varrho
(ac_k),
\end{equation}
where $a\in \kr(\e)$.

\smallskip
\bla{ii} Similarly, if $\Delta$ is an  (even)  derivation  on
$\hor_P$
then there exists the unique linear map
$\varrho\colon \cal{A}\rightarrow \hor_P$  such that
\begin{equation}\label{strD2}
\Delta(\varphi)=\sum_k\varphi_k\varrho(c_k).
\end{equation}
The following identities hold
\begin{align}
\Delta\varrho&=0\label{Delr}\\
\varrho(a)\varphi&=\sum_k\varphi_k\varrho(ac_k)\label{com2}
\end{align}
where $a\in\kr(\e)$.

\smallskip
\bla{iii} In both cases
\begin{align}
\varrho(1)&=0\label{r1}\\
(\varrho\otimes \id)\ad&=F^{\sstar}\varrho\label{rad}
\end{align}
where $\ad\colon \cal{A}\rightarrow\cal{A}\otimes\cal{A}$ is the adjoint
action of $G$ on itself. The map $\Delta$ is hermitian iff
\begin{equation}
\varrho(\kappa(a)^*)=-\varrho(a)^*\label{r*}
\end{equation}
for each $a\in\cal{A}$.
\end{lem}
\begin{pf}
We shall prove the statements assuming that $\Delta$ is
a derivation. The case when $\Delta$ is  an  antiderivation  can  be
treated similarly.

     As first, it is clear  that  $\Delta$  is  unique,  if  exists.
Indeed, property \eqref{strD2} implies
$$
\varrho(a)=\sum_i q_i\Delta(b_i),
$$
where $q_i,b_i\in\cal{B}$ are such that \eqref{free} holds.

     Conversely, the above formula can be taken  as  the  starting
point for a definition of $\varrho$. However, such a definition will be
consistent iff
$$\left\{\Sum_iq_iF(b_i)=0\right\}\,\Longrightarrow\,
\left\{\Sum_iq_i\Delta(b_i)=0\right\}.$$

     Here,  the  map  $\varrho$  will  be  constructed  in  a  slightly
different  way,  without  explicitly  proving   this   consistency
condition. The notation introduced in \cite{D2}--Appendix B will
be followed. Let $\cal T$ be the set  of  equivalence  classes  of
irreducible unitary representations of $G$.

For each class $\alpha\in\cal{T}$
let us consider an irreducible representation  $u^\alpha\in\alpha$
(acting in $\Bbb{C}^n$, with matrix elements $u^\alpha_{ij}$) and choose
elements    $b_{ki}^\alpha\in\cal{B}^\alpha$
(where $\cal{B}^\alpha$ is the     multiple
irreducible $\cal  V$-submodule  of  $\cal  B$,  corresponding  to
$\alpha$) as explained in Appendix B of
\cite{D2}
(in fact additional ``positivity'' assymptions made
there are not essential for the final result, in the general  case
only some minus signs should be added in appropriate places). Let
$\varrho\colon \cal{A}\rightarrow \hor_P$ be a linear map specified by
$\varrho(1)=0$ and
\begin{equation}\label{rdef}
\varrho(u_{ij}^\alpha)=\sum_kb^{\alpha*}_{ki}\Delta(b_{kj}^\alpha).
\end{equation}
Let  us  assume  that  elements $\xi_1,\dots,\xi_n\in \hor_P$
transform according to the  representation  $u^\alpha$.  In  other
words
$$F^\sstar(\xi_i)=\sum_j\xi_j\otimes u_{ji}^\alpha.$$
We have then
$$\sum_i\xi_i\varrho(u^\alpha_{ij})=\sum_{ki}\xi_ib_{ki}^{\alpha*}
\Delta
(b_{kj}^\alpha)=\sum_{ki}
\Delta(\xi_ib^{\alpha*}_{ki}b^{\alpha}_{kj}
)=\sum_i \delta_{ij}\Delta(\xi_i)=\Delta(\xi_j),
$$
because of \eqref{DelW} and the $F^{\sstar}$-invariance of
$\Sum_i\xi_ib^{\alpha*}_{ki}$. Consequently, property \eqref{strD2}
holds, because irreducible $\alpha$-multiplets span $\cal{B}^\alpha$,
and
$$\cal{B}=\sideset{}{^\oplus}\sum_{\alpha\in\cal{T}}
\cal{B}^\alpha.$$

     Let us check \eqref{rad}. Applying \eqref{FDel} and
\eqref{rdef} we obtain
\begin{multline*}
     F^{\sstar}\varrho(u^\alpha_{ij})=\sum_kF(b^{\alpha*}_{ki})
(\Delta\otimes \id)F(b^\alpha_{kj})\\
=\sum_{kmn}(b^{\alpha*}_{km}\otimes u^{\alpha*}_{mi})
\bigl(\Delta(b^\alpha_{kn})\otimes u^\alpha_{nj}\bigr)=\sum_{mn}
\varrho(u_{mn}^\alpha)
\otimes u_{mi}^{\alpha*}u_{nj}^{\alpha}
=(\varrho\otimes \id)\ad(u_{ij}^\alpha).
\end{multline*}
Acting by $(\Delta\otimes \id)$ on the identity
\begin{equation}\label{A6}
\sum_k b^{*\alpha}_{ki}F(b^\alpha_{kj})=1\otimes
u_{ij}^\alpha
\end{equation}
we find
$$
\sum_k\left(\Delta(b_{ki}^{\alpha*})F(b^\alpha_{kj})+b^{\alpha*}_{ki}(
\Delta\otimes
\id)F(b^{\alpha}_{kj})\right)=0.
$$
Identity \eqref{strD2} together with the above equality and \eqref{FDel}
gives
$$
     0=\sum_k\left(\Delta(b_{ki}^{\alpha*})\Delta(b^\alpha_{kj})+
b^{\alpha*}_{ki}\Delta^2(b^\alpha_{kj})\right)
=\Delta\biggl(\sum_kb^{\alpha*}_{ki}\Delta(b^\alpha_{kj})\biggr)=
\Delta\varrho(u_{ij}^\alpha).
$$
Hence, \eqref{Delr} holds. Multiplying \eqref{A6} by
$F^{\sstar}(\varphi)$  on
the right and using \eqref{strD2} we obtain
$$
\sum_{kml}b^{\alpha*}_{ki}b^\alpha_{km}\varphi_l\varrho
(u^\alpha_{mj}c_l)=\sum_kb^{\alpha*}_{ki}
\Delta(b^\alpha_{kj}\varphi)=
\sum_l\varphi_l\varrho(u^\alpha_{ij}c_l).
$$
On the other hand
$$
\varrho(u^\alpha_{ij})\varphi=\sum_kb^{\alpha*}_{ki}\Delta(
b^\alpha_{kj}  )\varphi=\sum_kb_{ki}^{\alpha*}\Delta(b^\alpha_{kj}
\varphi)-\delta_{ij}\Delta(\varphi).
$$
Consequently,
$$
\varrho(u_{ij}^\alpha)\varphi+\delta_{ij}\Delta(\varphi)
=\sum_l\varphi_l
\varrho(u_{ij}^\alpha c_l),
$$
which proves \eqref{com2}.

     Finally, let us assume that $\Delta$ is hermitian. This implies
\begin{multline*}
\varrho(u_{ij}^\alpha)^*=\biggl[\sum_kb_{ki}^{\alpha*}
\Delta(b^\alpha_{kj})\biggr]^*
=\sum_k \Delta(b^{\alpha*}_{kj})b^\alpha_{ki}\\
=-\sum_kb_{kj}^{\alpha*}\Delta(b^\alpha_{ki})=-\varrho(u_{ji}^\alpha)
=-\varrho(\kappa(u_{ij}^{\alpha})^*).
\end{multline*}
Hence, \eqref{r*} holds. Conversely, if \eqref{r*} holds
then  applying \eqref{strD2}, \eqref{com2} and \eqref{r1} we obtain
\begin{multline*}
\Delta(\varphi)^*=\sum_k\bigl(\varphi_k\varrho(c_k)\bigr)^*=
\sum_k\varrho(c_k)^*\varphi_k^*\\
=-\sum_k \varrho(\kappa(c_k)^*)\varphi_k^*=\sum_k\varphi_k^*\varrho(c_k^*)
=\Delta(\varphi^*)
\end{multline*}
which completes the proof.
\end{pf}

     It is worth noticing that degrees of $\varrho$ and $\Delta$
are the same. Further, $\Delta$
is  completely  fixed  by
its restriction on $\cal{B}$,
because $\varrho$ is expressible in terms of this restriction.

     It is possible to ``reverse'' the above construction of $\varrho$.
If a linear homogeneous map $\varrho\colon \cal{A}\rightarrow \hor_P$ is
given such that
\begin{equation}
\varrho(a)\varphi=(-1)^{\partial\varrho\partial\varphi}\sum_k
\varphi_k\varrho(ac_k)
\end{equation}
(where $a\in\kr(\e)$)
then \eqref{strD1} (or \eqref{strD2}) determines a map
$\Delta\colon \hor_P
\rightarrow \hor_P$,  which  is  an  even  (odd)  (anti)derivation,
depending on the parity of $\varrho$.

     Let us  now  assume  that $\Omega_M$
is  endowed  with  a differential *-algebra  structure,  specified
by a first-order differential map
$d_M\colon \Omega_M\rightarrow\Omega_M.$ In other words,  the  following
properties hold
\begin{gather}
d_M(\Omega_M^*)\subseteq\Omega_M^{*+1}\\
d_M(\varphi\psi)=d_M(\varphi)\psi+(-1)^{\partial\varphi}\varphi
d_M(\psi) \label{dMp}\\
d_M^2=0 \label{dM2}\\
d_M(\varphi^*)=d_M(\varphi)^*. \label{dM*}
\end{gather}
\begin{defn}
A {\it preconnection}  on  $P$  (with  respect   to
$\bigl\{\hor_P,F^\sstar,\Omega_M\bigr\}$) is a linear map
$D\colon \hor_P\rightarrow \hor_P$ satisfying
\begin{gather}
D(\hor_P^{*})\subseteq \hor_P^{*+1}\label{pre1}\\
F^{\sstar}D=(D\otimes \id)F^{\sstar}\label{pre3}\\
(D\restr\Omega_M)=d_M\label{pre2}\\
D(\varphi\psi)=D(\varphi)\psi+(-1)^{\partial\varphi}\varphi D(\psi)
\label{pre4}\\
D(\varphi^*)=D(\varphi)^*\label{pre5}.
\end{gather}
In  other  words,  preconnections  are  hermitian  right-covariant
first-order antiderivations on $\hor_P$, which extend $d_M$.
\end{defn}

     Let us observe that every linear map $D$  acting  in
$\hor_P$   and
satisfying \eqref{pre3} is reduced  in  the  space
$\Omega_M$.  If  a  map
$D\colon \hor_P\rightarrow \hor_P$ satisfies \eqref{pre1}--\eqref{pre4}
then $*D*$ possesses the same property, and hence
$(*D*+D)/2$ is a preconnection on $P$.

     Let us assume that $P$ admits preconnections. The set
$\pre(P)$ of all preconnections on $P$ is a  real  affine  space,  in
a  natural manner.  The  corresponding  vector  space
$\overrightarrow{\pre}(P)$
consists   of
hermitian first-order right-covariant antiderivations
$E$ on $\hor_P$ satisfying
\begin{equation}
 E(\Omega_M)=\{0\}.\label{EW}
\end{equation}
\begin{lem} \bla{i} For each $E\in\overrightarrow{\pre}(P)$
there exists the  unique
$\chi^\natural_E\colon \cal{A}\rightarrow \hor_P$
such that
\begin{equation}\label{strE}
E(\varphi)=-(-1)^{\partial\varphi}\sum_k \varphi_k
\chi^\natural_E(c_k), \end{equation}
for each $\varphi\in \hor_P$. We have
\begin{gather}
F^{\sstar}\chi^\natural_E(a)=(\chi^\natural_E\otimes \id)\ad(a)
\label{E1}\\
\chi^\natural_E(\kappa(a)^*)=-\chi^\natural_E(a)^*\label{E2}\\
\chi^\natural_E
(a)\varphi=(-1)^{\partial\varphi}\sum_k\varphi_k\chi^\natural_E
(ac_k)\label{E3}
\end{gather}
for each $a\in \kr(\e)$ and $\varphi\in \hor_P.$

\smallskip
\bla{ii} Similarly, for each $D\in\pre(P)$ there exists  the  unique
$\varrho^\natural_D\colon \cal{A}\rightarrow \hor_P$ satisfying
\begin{equation}\label{strDD}
D^2(\varphi)=-\sum_k\varphi_k\varrho^\natural_D(c_k),
\end{equation}
for each $\varphi\in \hor_P$. We have
\begin{gather}
F^{\sstar}\varrho^\natural_D(a)=(\varrho^\natural_D\otimes
\id)\ad(a)\label{DD1}\\
D\varrho^\natural_D(a)=0\label{DD2}\\
\varrho^\natural_D(\kappa(a)^*)=-\varrho^\natural_D(a)^*\label{DD3}\\
\varrho^\natural_D
(a)\varphi=\sum_k\varphi_k
\varrho^\natural_D
(ac_k)\label{DD4}
\end{gather}
for each $a\in\kr(\e)$.
\end{lem}
\begin{pf}  Let  us  consider  a  preconnection  $D$.   Properties
\eqref{pre1}--\eqref{pre5} imply that $D^2$  is  a
second-order right-covariant
hermitian derivation on $\hor_P$ satisfying
$$
D^2(\Omega_M)=\{0\}
$$
Applying Lemma~\ref{lem:21} ({\it ii}\/) to the case $\Delta=D^2$ we
conclude that there exists the unique map $\varrho^\natural_D
\colon \cal{A}\rightarrow \hor_P$ such that
\eqref{strDD} holds. Identities \eqref{DD1}--\eqref{DD4} follow from
statements ({\it ii-iii}\/)
in Lemma~\ref{lem:21}. The statement ({\it i}\/)
follows by a similar reasoning.
\end{pf}

We have $\chi^\natural_E(1)=\varrho^\natural_D(1)=0$.
Introduced maps are correlated such that
\begin{lem}\label{lem:23}
The following identity holds
\begin{equation}\label{r-nu}
\varrho^\natural_{D+E}(a)=
\varrho^\natural_D
(a)+D\chi^\natural_E(a)+
\chi^\natural_E(a^{(1)})\chi^\natural_E(a^{(2)}).
\end{equation}
\end{lem}
\begin{pf}
For  a  given  $a\in\cal{A}$  let  us  choose  elements
$q_i,b_i\in\cal{B}$ such that \eqref{free} holds,  and  let  us
assume  that
$F(b_i)=\Sum_kb_{ki}\otimes c_{ki}.$

A direct computation gives
\begin{equation*}
\begin{split}
-\varrho^\natural_{D+E}
(a)&=\sum_iq_i(D+E)^2(b_i)\\
&=-\sum_{ki}q_ib_{ki}\varrho^\natural_D(c_{ki})
-\sum_{ki}q_ib_{ki}\chi^\natural_E(c_{ki}^{(1)}
)\chi^\natural_E(c_{ki}^{(2)})\\
\phantom{=}&-\sum_{ki}q_iD\bigl(b_{ki}\chi^\natural_E(c_{ki})\bigr)
+\sum_{ki}q_i(Db_{ki})\chi^\natural_E(c_{ki})\\
&=-\varrho^\natural_D(a)-\chi^\natural_E(a^{(1)})
\chi^\natural_E(a^{(2)}
)-D\chi^\natural_E(a).
\qed
\end{split}
\end{equation*}
\renewcommand{\qed}{}
\end{pf}

     For each $D\in\pre(P)$ and $E\in\overrightarrow{\pre}(P)$
let $\cal{R}_D,\cal{P}_E\subseteq\kr(\e)$ be subspaces consisting of
elements anihilated by $\varrho^\natural_D$ and
$\chi^\natural_E$ respectively.
\begin{lem}\label{lem:24}
The spaces $\cal{R}_D$  and $\cal{P}_E$  are right
$\cal{A}$-ideals. Moreover,
\begin{gather}
\ad(\cal{P}_E)\subseteq \cal{P}_E\otimes\cal{A}\label{adRE}\\
\kappa(\cal{P}_E)^*=\cal{P}_E\label{kRE}\\
\kappa(\cal{R}_D)^*=\cal{R}_D\label{kRD}\\
\ad(\cal{R}_D)\subseteq\cal{R}_D\otimes\cal{A}.\label{adRD}
\end{gather}
\end{lem}
\begin{pf}
For  a  given $a\in\cal{A}$   let   us   choose   elements
$b_i,q_i\in\cal{B}$  such  that   \eqref{free}   holds.   Applying
\eqref{com2} we find
$$
\sum_iq_i\varrho^\natural_D(b)b_i=\sum_{ki}q_ib_{ki}
\varrho^\natural_D(bc_{ki})=\varrho^\natural_D(ba),
$$
for   each $b\in \kr(\e)$,   where $\Sum_k b_{ki}\otimes c_{ki}
=F(b_i).$ In
particular, if $b\in\cal{R}_D$  then $ba\in\cal{R}_D$, too.
Similarly, it follows that $\cal{P}_E$  are right
$\cal{A}$-ideals. Finally, \eqref{adRE}--\eqref{adRD}
directly follow from properties \eqref{E1}, \eqref{E2},
\eqref{DD1} and \eqref{DD3}.
\end{pf}
\section{Constructions of Differential Structures}

     In this section  two  constructions  will  be  presented.  As
first, starting from the system  of  preconnections   on   $P$  we
shall
construct a canonical  differential  calculus  on  the
structure quantum group $G$.  Secondly,  a  canonical  differential
structure on the bundle $P$ will be  constructed,  by  combining
this calculus  on  $G$  with  the  algebra
$\hor_P$.
As we  shall  see,  there
exists a natural correspondence between preconnections and regular
connections on $P$. Preconnections are  interpretable  as  covariant
derivatives associated to regular connections.

     Let ${\nat}$ be the intersection of all ideals
$\cal{R}_D$ and $\cal{P}_E$.
According to Lemma~\ref{lem:23}
\begin{equation}
{\nat}=\cal{R}_D{\Bcap}\biggl\{{\Bcap_E}\cal{P}_E
\biggr\},
\label{R2}
\end{equation}
for an arbitrary $D\in\pre(P)$.  Indeed,
\begin{equation}\label{prop-nu}
\chi^\natural_E(a^{(1)})\chi^\natural_E(a^{(2)})=\frac{1}{2}
\chi^\natural_E(a^{(2)})\chi^\natural_E(\kappa(a^{(1)})a^{(3)}),
\label{nua}
\end{equation}
as follows from \eqref{E1} and \eqref{E3}. In particular
if  $a$  belongs  to  the
right-hand side of \eqref{R2} then $\varrho^\natural_{D+E}(a)=0$,
for each $E\in\overrightarrow{\pre}(P)$.

     Let $\Psi$  be  the  left-covariant  first-order  differential
calculus on $G$ canonically corresponding to ${\nat}$
(in the sense of \cite{W-diff}). Lemma~\ref{lem:24} implies
\begin{gather}
\ad({\nat})\subseteq {\nat}
\otimes \cal{A}\label{bicR}\\
\kappa({\nat})^*= \nat. \label{*covR}
\end{gather}

In other words \cite{W-diff}, $\Psi$ is a bicovariant *-calculus.
     We shall assume that the complete  differential  calculus  on
$G$
is   based   on   the   universal    envelope    $\Psi^{\wedge}$ of
$\Psi$
(\cite{D1}--Appendix B).

     Let   us   consider    a    preconnection $D$,    and    let
$\varrho_D\colon \Psi_{\inv}\rightarrow \hor_P$  be a linear
map defined by
\begin{equation}
 \varrho_D\pi=\varrho^\natural_D, \label{rpi}
\end{equation}
where $\pi\colon \cal{A}\rightarrow \Psi_{\inv}$ is the canonical
projection map.
\begin{defn}
The map $\varrho_D$  is called {\it the curvature} of $D$.
\end{defn}
Similarly, for each $E\in\overrightarrow{\pre}(P)$ let
$\chi_E\colon \Psi_{\inv}\rightarrow \hor_P$ be a map given by
\begin{equation}
\chi_E\pi=\chi^\natural_E.\label{nupi}
\end{equation}

     The following identities summarize results of the previous
section:
\begin{gather*}
\varrho_{D+E}=\varrho_D+D\chi_E-\langle
\chi_E,\chi_E\rangle\qquad\quad D\varrho_D=0\\
F^{\sstar}\chi_E=(\chi_E\otimes \id)\Ad\qquad\quad
\chi_E(\vartheta^*)=\chi_E(\vartheta)^*\\
F^{\sstar}\varrho_D=(\varrho_D\otimes \id)\Ad\qquad\quad
\varrho_D(\vartheta^*)=\varrho_D(\vartheta)^*\\
\begin{aligned}
-D^2(\varphi)=&\sum_k\varphi_k\varrho_D\pi(c_k)\\
\varrho_D(\vartheta)\varphi=&\sum_k\varphi_k\varrho_D(\vartheta
{\circ} c_k)
\end{aligned}\qquad
\begin{aligned}
-E(\varphi)=&(-1)^{\partial\varphi}\sum_k\varphi_k\chi_E\pi(c_k)\\
\chi_E(\vartheta)\varphi=&(-1)^{\partial\varphi}
\sum_k\varphi_k\chi_E(\vartheta{\circ} c_k)
\end{aligned}
\end{gather*}
where $\langle,\rangle$ are the brackets
associated  to  an  arbitrary  ``embedded
differential'' map $\delta\colon\Psi_{\inv}\rightarrow\Psi_{\inv}
\otimes\Psi_{\inv}$, and $\varpi\colon\Psi_{\inv}\rightarrow
\Psi_{\inv}\otimes\cal{A}$ is the (co)adjoint action of $G$ on
$\Psi_{\inv}$. The $\circ$ denotes the natural right $\cal{A}$-module
structure on $\Psi_{\inv}$.

Let us now consider a *-algebra $\vh_P$ representing
``vertically-horizontally'' decomposed forms
\cite{D2}. At the level of (graded) vector
spaces
$$
\vh_P=\hor_P\otimes \Psi_{\inv}^\wedge,
$$
where $\Psi_{\inv}^\wedge\subseteq\Psi^\wedge$
is the subalgebra of left-invariant elements.
     The *-algebra structure on $\vh_P$  is specified by
\begin{align}
(\varphi\otimes\vartheta)^* &=\sum_k
\varphi_k^*\otimes(\vartheta^* {\circ} c_k^*)\label{*-str1}\\
(\psi\otimes\eta)(\varphi\otimes\vartheta)&=
(-1)^{\partial\eta\partial\varphi}\sum_k
\psi\varphi_k\otimes(\eta{\circ} c_k)\vartheta.\label{*-str2}
\end{align}
By construction
$\hor_P$ and $\Psi_{\inv}^\wedge$ are interpretable
as subalgebras of $\vh_P$.

The formulas
\begin{align}
\partial_D(\varphi)&=D(\varphi)+(-1)^{\partial\varphi}
\sum_k\varphi_k \pi(c_k)\label{for1}\\
\partial_D(\vartheta)&=\varrho_D(\vartheta)+d(\vartheta)\label{for2}
\end{align}
(where $\varphi\in\hor_P$ and $\vartheta\in\Psi_{\inv}$, while
$d\colon \Psi_{\inv}^\wedge
\rightarrow\Psi_{\inv}^\wedge$ is the corresponding
differential) determine
(via the graded Leibniz rule) a hermitian first-order differential
$\partial_D\colon\vh_P\rightarrow\vh_P$, for each $D\in\pre(P)$.

Differential   *-algebras
$(\vh_P,\partial_D)$ are
mutually naturally isomorphic.
\begin{pro}\label{pro:iso} \bla{i} For each $E\in\overrightarrow{\pre}(P)$
there  exists  the
unique homomorphism
$h_E\colon \vh_P\rightarrow\vh_P$ such that
\begin{align}
h_E(\varphi)&=\varphi\label{hE1}\\
h_E(\vartheta)&=\vartheta-\chi_E(\vartheta), \label{hE2}
\end{align}
for each $\varphi\in \hor_P$ and $\vartheta\in\Psi_{\inv}$.

\smallskip
\bla{ii} The maps $h_E$ are hermitian and bijective.

\smallskip
 \bla{iii} We have
\begin{equation}
h_Eh_W=h_{E+W}\label{hE4}
\end{equation}
for each $E,W\in\overrightarrow{\pre}(P).$

\smallskip
\bla{iv} The  map $h_E$   is  an  isomorphism  between
differential structures $(\vh_P,\partial_D)$ and
$(\vh_P,\partial_{D+E})$,
for  each $E\in\overrightarrow{\pre}(P)$ and
$D\in\pre(P)$.
\end{pro}
\begin{pf}
Uniqueness of maps $h_E$  follows from the  fact  that
$\hor_P$ and $\Psi_{\inv}$ generate $\vh_P$.
To  establish their existence,  it  is
sufficient to check that conditions \eqref{hE1}--\eqref{hE2} are
compatible  with
the product rule  for
$\vartheta\varphi$,  and  with  the  quadratic
constraint defining the algebra $\Psi_{\inv}^{\wedge}$. We have
\begin{multline*}
\vartheta\varphi=(-1)^{\partial\varphi}\sum_k\varphi_k
(\vartheta{\circ} c_k)\longrightarrow
(-1)^{\partial\varphi}\sum_k \bigl(\varphi_k
(\vartheta{\circ} c_k)-
\varphi_k \chi_E (\vartheta{\circ} c_k)\bigr)\\
=\bigl[\vartheta-\chi_E(\vartheta)\bigr]\varphi,
\end{multline*}
for each $\vartheta\in\Psi_{\inv}$ and $\varphi\in
\hor_P$.
Further, if $a\in\nat$ then
\begin{multline*}
0=\pi(a^{(1)})\pi(a^{(2)}) \longrightarrow
\pi(a^{(1)})\pi(a^{(2)})+\chi_E\pi(a^{(1)})\chi_E\pi(a^{(2)})\\
-\bigl[\chi_E\pi(a^{(1)})\bigr]\pi(a^{(2)})-
\pi(a^{(1)})\chi_E\pi(a^{(2)})=0,
\end{multline*}
because of \eqref{prop-nu} and
\begin{equation*}
\begin{split}
\pi(a^{(1)})\chi_E\pi(a^{(2)})&=-\chi_E\pi(a^{(3)})
\bigl[\pi(a^{(1)}
){\circ} \bigl(\kappa(a^{(2)})a^{(4)}\bigr)\bigr]\\
&=-\bigl[\chi_E\pi(a^{(1)})\bigr]\pi(a^{(2)})+
\chi_E\pi(a^{(2)})\pi\bigl[\kappa(a^{(1)})a^{(3)}\bigr]\\
&=-\bigl[\chi_E\pi(a^{(1)})\bigr]\pi(a^{(2)}).
\end{split}
\end{equation*}
Hence, $h_E$  exists.
     In order to prove the hermicity of $h_E$  it  is  sufficient  to
check that restrictions of $h_E$ on $\Psi_{\inv}$ and $\hor_P$ are
hermitian maps.
This imediately follows from the hermicity of $\chi_E$.
Similarly, it  is  sufficient
to check that \eqref{hE4} holds on $\hor_P$ and $\Psi_{\inv}$.
This  trivially
follows from \eqref{hE1}--\eqref{hE2}. Now \eqref{hE4} implies that
$h_E$ are bijective maps.

     Let us prove ({\it iv}\/). Because of the graded Leibniz rule,  it
is sufficient to check that $h_E\partial_D =\partial_{D+E}h_E$ holds on
$\Psi_{\inv}$ and $\hor_P$. We have
\begin{equation*}
\begin{split}
h_E\partial_D(\varphi)&=D(\varphi)+(-1)^{\partial\varphi}\sum_k
\varphi_k\bigl[\pi(c_k)-\chi_E\pi(c_k)\bigr]\\
&=(D+E)(\varphi)+(-1)^{\partial\varphi}
\sum_k\varphi_k\pi(c_k)=\partial_{D+E}(\varphi).
\end{split}
\end{equation*}
Further,
\begin{equation*}
\begin{split}
h_E\partial_D\pi(a)&=-h_E\bigl[\pi(a^{(1)})\pi(a^{(2)})\bigr]+
\varrho_D\pi(a)\\
&=-\pi(a^{(1)})\pi(a^{(2)})-\chi_E\pi(a^{(1)})\chi_E\pi(a^{(2)})\\
&\phantom{=}+\bigl[\chi_E\pi(a^{(1)})\bigr]\pi(a^{(2)})+
\pi(a^{(1)})\chi_E\pi(a^{(2)})+\varrho_D\pi(a)\\
&=-\pi(a^{(1)})\pi(a^{(2)})-\chi_E\pi(a^{(1)})\chi_E\pi(a^{(2)})\\
&\phantom{=}+\chi_E \pi(a^{(2)})\pi\bigl[\kappa(a^{(1)})a^{(3)}\bigr]+
\varrho_D\pi(a)\\
&=\varrho_D\pi(a)+D\chi_E\pi(a)+\chi_E\pi(a^{(1)})\chi_E\pi(a^{(2)}
)\\
&\phantom{=}-D\chi_E\pi(a)-\pi
(a^{(1)})\pi(a^{(2)})\\
&\phantom{=}+\chi_E\pi(a^{(2)})\pi\bigl[\kappa(a^{(1)})a^{(3)}\bigr]-
\chi_E\pi(a^{(2)})\chi_E\pi\bigl[\kappa(a^{(1)})a^{(3)}\bigr]\\
&=\varrho_{D+E}\pi(a)-D\chi_E\pi(a)-
\pi(a^{(1)})\pi(a^{(2)})\\
&\phantom{=}+\chi_E\pi(a^{(2)})\pi\bigl[\kappa(a^{(1)})a^{(3)}\bigr]
-\chi_E\pi(a^{(2)})\chi_E\pi\bigl[\kappa(a^{(1)})a^{(3)}\bigr]\\
&=\partial_{D+E}\pi(a)-D\chi_E\pi(a)-E\chi_E\pi(a)
+\chi_E\pi(a^{(2)})\pi\bigl[\kappa(a^{(1)})a^{(3)}\bigr]\\
&=\partial_{D+E}\pi(a)-\partial_{D+E}\chi_E\pi(a)=
\partial_{D+E}h_E\pi(a).\qed
\end{split}
\end{equation*}
\renewcommand{\qed}{}
\end{pf}

     Now, a {\it manifestly invariant}
differential calculus on $P$ can be constructed by ``gluing''
algebras $(\vh_P,\partial_D)$, with the help of isomorphisms $h_E$. Let
$\Omega_P$ be a graded-differential *-algebra obtained in this way.
By construction, each $D\in\pre(P)$ naturally determines a *-isomorphism
$\pi_D\colon\Omega_P\rightarrow\vh_P$, such that
$$ \pi_D d_P=\partial_D\pi_D$$
where $d_P$ is the differential on $\Omega_P$. We have
$$ h_E=\pi_{D+E}\pi_{D}^{-1}$$
for each $D\in\pre(P)$ and $E\in\overrightarrow{\pre}(P)$.

The map $F$ is naturally extendible to a homo\-mor\-phism
$\widehat{F}\colon\Omega_P\rightarrow
\Omega_P\grten
\Psi^\wedge$ of graded-differential *-algebras. Explicitly,
$$ (\pi_D^{\phantom{1}}\widehat{F}\pi_D^{-1})(\varphi\otimes
\vartheta)=F^\sstar(\varphi)\widehat{\Ad}(\vartheta)$$
where $\widehat{\Ad}\colon
\Psi_{\inv}^\wedge\rightarrow
\Psi_{\inv}^\wedge\grten\Psi^\wedge$ is the
graded-differential *-homo\-mor\-phism extending the (co)adjoint action
$\Ad
\colon\Psi_{\inv}\rightarrow\Psi_{\inv}\otimes\cal{A}$.

The following equality holds
$$ \hor_P=(\widehat{F})^{-1}\Bigl\{\Omega_P\otimes
\cal{A}\Bigr\}.$$
This justifies the interpretation of $\hor_P$, as the algebra of
horizontal forms. Also, the above equality justifies the interpretation of
$\Omega_M$, as consisting of differential forms on the base manifold $M$. In
particular, $\Omega_M$ is $d_P$-invariant, and $(d_P\restr\Omega_M)=d_M$.

\begin{pro}
For each $D\in\pre(P)$, the connection $\omega=\omega_D$ given by
\begin{equation}
\omega(\vartheta)=\pi_D^{-1}(1\otimes\vartheta)\label{w-D}
\end{equation}
is regular (and multiplicative). Moreover,
\begin{gather}
D=D_{\omega}\\
\pi_D^{-1}=m_{\omega}\label{pi-dec}
\end{gather}
where $m_{\omega}\colon \hor_P\otimes\Psi_{\inv}^{\wedge}
\rightarrow\Omega_P$ is the
canonical decomposition map.

Conversely, if $\omega$ is an arbitrary
regular connection on $P$ then
there exists the unique $D\in\pre(P)$ such that \eqref{w-D} holds.
\end{pro}
\begin{pf} The first part of the proposition directly follows from
the definition of $\omega$.
Let
$\omega\colon\Psi_{\inv}\rightarrow\Omega_P$  be  an
arbitrary regular connection.
Let $D\colon \hor_P\rightarrow \hor_P$ be the covariant  derivative
associated to $\omega$.  By  construction,  $D\in\pre(P)$  and
\eqref{pi-dec} holds. This further implies that
\eqref{w-D} holds. On the
other hand, if \eqref{w-D} holds then $D=D_{\omega}.$
\end{pf}

     In other words, formula \eqref{w-D} establishes
a bijective affine correspondence between the spaces of
regular connections and preconnections.

\section{Concluding Remarks}

     It is worth noticing that the presented construction  of  the
calculus  on  the  bundle  works for  an   arbitrary
bicovariant  first-order  *-calculus $\Gamma$ based  on   a   right
$\cal{A}$-ideal $\cal{R}$ satisfying
$\cal{R}\subseteq\nat$. It is also possible
to perform the construction of differential structures on $G$ and $P$,
dealing with a restricted set of preconnections, forming an appropriate
affine subspace $\cal{L}\subseteq\pre(P)$.
Covariant derivatives of regular
connections (with respect to the associated
differential structures) form an affine subspace
$\cal{L}^\sstar\supseteq\cal{L}$ of $\pre(P)$. Particularly
interesting are subspaces $\cal{L}$ satisfying the
``stability property'' $\cal{L}=\cal{L}^\sstar$.

     Let us turn back  to  the  general  context  of  differential
calculus on quantum principal bundles. Let us assume that the calculus
on $G$ is based on (the universal envelope of) a first-order bicovariant
*-calculus $\Gamma$ (determined by a right $\cal{A}$-ideal $\cal{R}$).

Let $P=(\cal{B},i,F)$  be  a
quantum principal bundle,
and let  us  assume  that  the  calculus  on $P$   is   based
on   a
graded-differential *-algebra $\Omega(P)$. Further,  let  us  assume
that $P$ admits regular multiplicative connections
(relative to $\Omega(P)$).

Then for every regular connection $\omega$ the corresponding
covariant derivative map $D_{\omega}\colon \hor(P)\rightarrow
\hor(P)$
is   a   preconnection   on   $P$   (relative   to
$\bigl\{\hor(P),F^{\wedge},\Omega(M)\bigr\}$).
However, the converse is generally not  true.
In general only an affine subspace of $\pre(P)$ will be  induced  by
regular connections. On the other hand, starting from  $\hor(P)$  and
$\pre(P)$ and  applying  constructions  presented  in  this
study we obtain the bicovariant *-calculus $\Psi$ on
$G$ and the
graded-differential  *-algebra $\Omega_P$.   In
general   algebras
$\Omega(P)$ and $\Omega_P$  are not mutually naturally
related. However,
if all preconnections are interpretable as  covariant  derivatives
of regular connections (relative  to $\Omega(P)$)  then  (and  only
then) $\cal{R}\subseteq\nat$. In this case
$\Omega_P$ (and $\Psi^{\wedge}$)
can    be    obtained    by    factorizing     $\Omega(P)$ (and
$\Gamma^\wedge$)   through
appropriate graded-differential *-ideals.

     Presented constructions of differential structures generalize
the    corresponding    constructions      of      the      theory
\cite{D1} of  quantum principal bundles over classical smooth
manifolds.

     Let us assume that $P$ is a quantum principal $G$-bundle  over  a
(compact) smooth manifold $M$ (in the sense of  \cite{D1}).  The  algebra
$\hor(P)$  representing  horizontal  forms  can  be  constructed   by
combining  differential  forms  on  $M$    with    functions    on
$P$
(indenpendently of the specification of the complete calculus on the
bundle and the structure group).
If the bundle is locally trivialized over some open  set,
then $\hor(P)$ will be locally trivialized in a natural manner, too.

     It turns out that there exists
a  natural bijection between preconnections on $P$, and standard
connections  on
the  classical part $P_{cl}$ of $P$.
     The right $\cal{A}$-ideal $\nat$
consists   precisely   of   those
elements $a\in \kr(\e)$  satisfying
$$(X\otimes \id)\ad(a)=0,$$
  for   each
$X\in \lie(G_{cl})$ (where the elements of $\lie(G_{cl})$
are  understood  as
(hermitian) functionals $X\colon \cal{A}\rightarrow \Bbb{C}$
satisfying
$X(ab)=\e(a)X(b)+X(a)\e(b)$).

     Hence, $\nat$ determines the minimal admissible
(bicovariant *-) calculus on $G$ (in
the terminology of \cite{D1}). In other words, $\Psi$ is
the   minimal   left-covariant   first-order   calculus   on   $G$
compatible,    in   a   natural    manner,    with    all    local
retrivializations of the bundle $P$. Let
$\Omega(P)$ be the graded-differential  *-algebra  constructed  from
$\Psi$ and $P$, with the help of
$G$-cocycles \cite{D1}. Then the identity map on $\cal B$  extends
to the graded-differential *-isomorphism between $\Omega_P$
and $\Omega(P)$.

\appendix
\section{Differential Structures Based On Exterior Algebras}
Let $\sigma\colon\Psi{\otimes_{\cal A}}
\Psi\rightarrow\Psi{\otimes_{\cal A}}\Psi$ be the canonical
flip-over operator \cite{W-diff}. This map is
a bicovariant bimodule automorphism. Its ``left-invariant'' part
$\sigma\colon\Psi_{\inv}^{\otimes 2}\rightarrow
\Psi_{\inv}^{\otimes 2}$ is given by
\begin{equation}\label{flip}
\sigma(\eta\otimes\vartheta)=\sum_k\vartheta_k\otimes
(\eta{\circ} c_k)
\end{equation}
where $\Sum_k \vartheta_k\otimes c_k=\Ad(\vartheta)$.
Let $\Psi^\vee$ be the corresponding
exterior algebra \cite{W-diff}.
By definition, $\Psi^\vee$
can be obtained by factorizing the ``tensor bundle'' algebra
$\Psi^\otimes$ through the (bicovariant *-) ideal $S^\vee=\kr(A)$
where
$$A={\sum_{n\geq 0}}^\oplus A_n$$
is the corresponding ``total antisymetrizer''.
The maps $A_n\colon\Psi^{\otimes n}\rightarrow\Psi^{\otimes n}$
are given by
$$ A_n=\sum_{\pi\in S_n}(-1)^\pi \sigma_\pi.$$
Here $\sigma_\pi\colon\Psi^{\otimes n}\rightarrow
\Psi^{\otimes n}$ are obtained by replacing
transpositions
in a minimal decomposition of $\pi$ with the corresponding
$\sigma$-twists.

The space $\hor_P\otimes\Psi_{\inv}^\vee=\vh_P^\vee$ possesses
a natural *-algebra structure (expressed by the same formulas
\eqref{*-str1}--\eqref{*-str2}). As shown in \cite{D2}--Appendix A
the formulas \eqref{for1}--\eqref{for2} consistently determine
a first-order hermitian differential
$\partial_D$ on $\vh_P^\vee$, for each $D\in\pre(P)$.
Now we shall
prove that (the analog of Proposition~\ref{pro:iso}) differential
algebras $(\vh_P^\vee,\partial_D)$ are naturally isomorphic.

Let us consider the *-algebra
$\amalg_P=\hor_P\otimes\Psi_{\inv}^\otimes$ (the *-structure is specified
by   the same formulas as for $\vh_P$).   For   each
$E\in\overrightarrow{\pre}(P)$ there exists the unique *-automorphism
$h_E^\sstar\colon \amalg_P\rightarrow \amalg_P$
satisfying
\begin{align*}
h^\sstar_E(\varphi)&=\varphi\\
h^\sstar_E(\vartheta)&=\vartheta-\chi_E(\vartheta)
\end{align*}
for each $\varphi\in\hor_P$ and $\vartheta\in\Psi_{\inv}$.
Moreover
$$
h^\sstar_{E+W}=h^\sstar_Eh^\sstar_W
$$
for each $E,W\in\overrightarrow{\pre}(P)$.
We shall prove that
$$h^\sstar_E(\hor_P\otimes S^\vee_{\inv}
)=\hor_P\otimes S^\vee_{\inv}$$
for each $E\in\overrightarrow{\pre}(P)$. Applying \eqref{flip} and
elementary properties of maps $\chi_E$ we find
\begin{equation}
h_{-E}^\sstar(\vartheta)=\sum_{k+l=n}(\chi_E^{\otimes k}\otimes \id^l)
A_{kl}(\vartheta)=\sum_{k+l=n}\frac{1}{k!}(\chi_E^{\otimes k}A_k\otimes
\id^{l})A_{kl}(\vartheta)
\end{equation}
for each $\vartheta\in\Psi_{\inv}^{\otimes n}$,
where $\chi_E^{\otimes}\colon\Psi_{\inv}^\otimes\rightarrow
\hor_P$ is  the  corresponding  unital  multiplicative
extension, and
$$ A_{kl}=\sum_{\pi\in S_{kl}} (-1)^\pi \sigma_{\pi^{-1}}.$$
Here $S_{kl}\subseteq  S_{k+l}$  is  consisting   of   permutations
preserving the order of sets $\{1,\dots,k\}$
and $\{k+1,\dots,k+l\}$,
(and $\sigma_\pi$ are restricted in $\Psi_{\inv}^\otimes$).

In particular if $\vartheta\in S^\vee_{\inv}$ then $h_E^\sstar(\vartheta)
\in\hor_P\otimes S^\vee_{\inv}$, because of
$$ A_{k+l}=(A_k\otimes A_{l})A_{kl}. $$
 Thus,  the  maps  $h^\sstar_E$  can  be
``factorized'' through $S^\vee_{\inv}$.  In  such  a  way  we  obtain
*-automorphisms
$h_E\colon\vh_P^\vee\rightarrow\vh_P^\vee$. We have
\begin{equation*}
h_E \partial_D=\partial_{D+E}h_E
\end{equation*}
for each $D\in\pre(P)$ (for simplicity, we have denoted by the same symbols
basic maps operating in $\vh_P$ and $\vh^\vee_P$).
This follows from Proposition~\ref{pro:iso}
({\it iv}\/), and from the fact that there exists a  natural  differential
algebra epimomorphism
$\Uni\colon\Psi^\wedge\rightarrow\Psi^\vee$
reducing to  the
identity on $\cal A$ and $\Psi$ (a consequence  of the
universality of $\Psi^\wedge$). The algebra $\vh^\vee_P$
can be obtained
by factorizing $\vh_P$ through a graded *-ideal
$$\I_P=
\hor_P\otimes[S_{\inv}^\vee]^\wedge.$$
The rest of the construction
of the corresponding invariant
calculus on the bundle is the same as in the universal case.
In such a way  we  obtain  a  graded-differential  *-algebra
$\Lambda_P$. We have
$$\Lambda_P=\Omega_P/\I_P^\sstar$$
where $\I_P^\sstar$ is a  graded-differential *-ideal corresponding  to
$\I_P$ (the ideal $\I_P$ is $h_E$-invariant).
The algebra $\Lambda_P$ can also be constructed
directly from $\Omega_P$, using a general construction  described
in \cite{D2}--Appendix A.


\begin{thebibliography}{4}
\bibitem{D1} {\Dj}ur{\dj}evi\'c M:
{\it Geometry of Quantum Principal Bundles I}, Preprint QmmP 6/92,
Belgrade University, Serbia
\bibitem{D2} {\Dj}ur{\dj}evi\'c M:
{\it Geometry of Quantum Principal Bundles II--extended version},
Preprint, Instituto de Matematicas, UNAM, M\'exico (1994)
\bibitem{W-cmpg} Woronowicz S L: {\it Compact matrix  pseudogroups}
CMP {\bf 111} 613--665 (1987)
\bibitem{W-diff} Woronowicz S L: {\it Differential calculus on
compact matrix pseudogroups (quantum groups)} CMP {\bf 122}
125--170 (1989)
\end{thebibliography}
\end{document}